\documentclass{amsart}
\usepackage{amssymb,latexsym}
\numberwithin{equation}{section}
\newtheorem{theorem}{Theorem}[section]

\newtheorem{corollary}[theorem]{Corollary}

\newtheorem{remark}[theorem]{Remark}

\newtheorem{example}[theorem]{Example}

\newtheorem{notation}[theorem]{Notation} 
\begin{document}

\pagenumbering{arabic}
\pagestyle{headings}
\def\sof{\hfill\rule{2mm}{2mm}}
\def\ls{\leq}
\def\gs{\geq}
\def\SS{\mathcal S}
\def\qq{{\bold q}}
\def\txx{{\frac1{2\sqrt{x}}}}
\def\mn{\mbox{-}}

\title{ {\sc restricted $1\mn3\mn2$ permutations and generalized patterns}}
   
\author{Toufik Mansour}
\maketitle
\begin{center}{LaBRi, Universit\'e Bordeaux I,\\ 
              351 cours de la Lib\'eration, 33405 Talence Cedex\\
		{\tt toufik@labri.fr} }
\end{center}
%
\section*{Abstract}
Recently, Babson and Steingrimsson (see \cite{BS}) introduced generalized permutations 
patterns that allow the requirement that two adjacent letters in a 
pattern must be adjacent in the permutation. 

We study generating 
functions for the number of permutations on $n$ letters avoiding 
$1\mn3\mn2$ (or containing $1\mn3\mn2$ exactly once) 
and an arbitrary generalized pattern $\tau$ on $k$ letters, 
or containing $\tau$ exactly once. In several cases the generating 
function depends only on $k$ and is expressed via Chebyshev polynomials 
of the second kind, and generating function of Motzkin numbers.
\section{Introduction}
{\bf Permutation patterns:} let $\alpha\in S_n$ and $\tau\in S_k$ be two permutations. We say that $\alpha$ 
{\em contains} $\tau$ if there exists a subsequence $1\leq i_1<i_2<\cdots<i_k\leq n$ 
such that $(\alpha_{i_1},\dots,\alpha_{i_k})$ is order-isomorphic to $\tau$; 
in such a context $\tau$ is usually called a {\em pattern}. We say $\alpha$ {\em avoids}
$\tau$, or is $\tau$-{\em avoiding}, if such a subsequence does not exist. 
The set of all $\tau$-avoiding permutations in $S_n$ is denoted $S_n(\tau)$. 
For an arbitrary finite collection of patterns $T$, we say that $\alpha$
avoids $T$ if $\alpha$ avoids any $\tau\in T$; the corresponding subset of 
$S_n$ is denoted $S_n(T)$.

While the case of permutations avoiding a single pattern has attracted
much attention, the case of multiple pattern avoidance remains less
investigated. In particular, it is natural, as the next step, to consider
permutations avoiding pairs of patterns $\tau_1$, $\tau_2$. This problem
was solved completely for $\tau_1,\tau_2\in S_3$ (see \cite{SS}), for 
$\tau_1\in S_3$ and $\tau_2\in S_4$ (see \cite{W}), and for 
$\tau_1,\tau_2\in S_4$ (see \cite{B1,Km} and references therein).   
Several recent papers \cite{CW,MV1,Kr,MV3,MV2} deal with the case 
$\tau_1\in S_3$, $\tau_2\in S_k$ for various pairs $\tau_1,\tau_2$. Another
natural question is to study permutations avoiding $\tau_1$ and containing
$\tau_2$ exactly $t$ times. Such a problem for certain $\tau_1,\tau_2\in S_3$ 
and $t=1$ was investigated in \cite{R}, and for certain $\tau_1\in S_3$, 
$\tau_2\in S_k$ in \cite{RWZ,MV1,Kr,MV3}. The tools involved in these papers 
include continued fractions, Chebyshev polynomials, and Dyck paths.\\

{\bf Generalized permutation patterns:}
in \cite{BS} Babson and Steingrimsson introduced generalized permutation 
patterns that allow the requirement that two adjacent letters in a pattern
must be adjacent in the permutation. The idea for Babson and Steingrimsson in 
introducing these patterns was study of Mahonian statistics, and they proved that 
all Mahonian permutation statistics in the literature can be written as linear 
combinations of such patterns. 

We define our {\em generalized patterns} as words in the alphabet $1,2,3,\dots$, where 
two adjacent letters may or may not be separated by a dash. The absence of a dash 
between two adjacent letters in a pattern indicates that the corresponding letters 
in the permutation must be adjacent, and in the order (order-isomorphic) 
given by the pattern.
An example of a generalized pattern is $1\mn23$. 
An $1\mn23$ subword of a permutation $\pi=(\pi_1,\pi_2,\cdots,\pi_n)$ is a subword 
$(\pi_i,\pi_j,\pi_{j+1})$ ($i<j$), such that $\pi_i<\pi_j<\pi_{j+1}$.

Claesson \cite{C} presented a complete solution for the number of permutations avoiding any 
single generalized pattern of length three with exactly one adjacent pair of letters, 
and presented some results for the number of permutations avoiding two patterns 
of length three with exactly one adjacent pair of letters. Claesson and Mansour \cite{CM}
presented a complete solution for the number of permutations avoiding any double 
generalized patterns of length three with exactly one adjacent pair of letters. 
Besides, Kitaev \cite{Ki} investigate simultaneous avoidance of two or more $3$-letters 
generalized patterns without internal dashes.\\

A {\em Motzkin path} of length $n$ is a lattice path from $(0,0)$ to $(n,0)$ 
with steps $(1,0)$, $(1,1)$, and $(1,-1)$ that never goes below the $a$-axis. 
Letting $L$, $U$, and $D$ represent the steps $(1,0)$, $(1,1)$, and $(1,-1)$ 
respectively, we code such a path with word over $\{L,U,D\}$. 
If $p$ any nonempty Motzkin path, then $p$ can be decomposed as either, 
$p=Lp'$ or $p=Up'Dq'$ where $p',q'$ are Motzkin paths. 
The $n$th Motzkin number $M_n$ is the number of Motzkin paths of length $n$. 
The first few of the $M_n$ are 
	$$1,1,2,4,9,21,51,127,323,835,2188,5798,15511,41835,113634,310572.$$

\begin{notation}
The generating function for Motzkin numbers $M_n$ we denote by $M(x)$; that is, 
$M(x)=\sum_{n\geq 0} M_nx^n$.
\end{notation}

A {\em Chebyshev polynomials of the second kind} are defined $U_p(\cos\theta)=\sin(p+1)\theta/\sin\theta$. 
Following \cite{CW,MV1,Kr,MV3,MV2}, introduce 
	$$R_p(x)=\frac{U_{p-1}\left( \frac{1}{2\sqrt{x}} \right) }{\sqrt{x} U_p\left( \frac{1}{2\sqrt{x}} \right) }.$$ 

As consequence to \cite{MV3}, 
we present a general approach to the study of permutations 
avoiding $1\mn3\mn2$ and avoiding a generalized pattern in $S_k$ (or containing 
exactly once). As a consequence, we derive many of the previously known results 
for this kind of problems, as well as many new results.

The paper is organized as follows. The case of permutations avoiding both $1\mn3\mn2$ 
and $\tau$ generalized pattern is treated in Section $2$. We derive a simple recursion for the 
corresponding generating function for general $\tau$. This recursion can be solved explicitly 
for several interesting cases. The case of permutations avoiding $1\mn3\mn2$ and containing 
$\tau$ once is treated in Section $3$. Here again we start from general recursion, and then 
solve it for several particular cases. The case containing $132$ 
exactly once is treate in section $5$ and $6$ for avoiding $\tau$, or containing 
$\tau$ exactly once; respectively.

Most of the explicit solutions obtained in Sections $2$ and $3$ involve Chebyshev 
polynomials of the second kind and generating function of Motzkin numbers.
\section{Avoiding $1\mn3\mn2$ and a generalized pattern}

Consider an arbitrary pattern $\phi=(\phi_1,\dots,\phi_k)\in S_k(1\mn3\mn2)$. 
Recall that $\phi_i$ is said to be a {\it right-to-left maximum} if 
$\phi_i>\phi_j$ for any $j>i$. Let $m_0=k,m_1,\dots,m_r$
be the right-to-left maxima of $\phi$ written from left to right. 
Then $\phi$ can be represented as
	$$\phi=(\phi^0,m_0,\phi^1,m_1,\dots,\phi^r,m_r),$$
where each of $\phi^i$ may be eventually empty, and all the entries of $\phi^i$
are greater than all the entries of $\phi^{i+1}$. This representation is called
the {\it canonical decomposition\/} of $\phi$. Given the canonical
decomposition, we define the $i$th {\it prefix\/} of $\phi$ by 
$\pi^i=(\phi^0,m_0,\dots,\phi^i,m_i)$ for $1\leq i\leq r$ and $\pi^0=\phi^0$,
$\pi^{-1}=\emptyset$.
Besides, the $i$th {\it suffix\/} of $\phi$ is defined by
$\sigma^i=(\phi^i,m_i,\dots,\phi^r,m_r)$ for $0\leq i\leq r$ and
$\sigma^{r+1}=\emptyset$. The following representation 
$\tau=\phi^0\mn m_0\mn\phi^1\mn m_1\mn\cdots\mn\phi^r\mn m_r$ is called the {\it canonical decomposition for generalized pattern} 
$\phi\in S_k(1\mn3\mn2)$. Similarly we define $\pi^i$ and $\sigma^i$ for $\tau$.

Let $f_\tau(n)$ denote the number of permutations in $S_n(1\mn3\mn2)$ such that avoiding 
$\tau$, and let $F_\tau(x)=\sum_{n\geq 0}f_\tau(n)x^n$ be the corresponding
generating function. By $f_\tau^\beta(n)$ we denote the number of permutations
in $S_n(1\mn3\mn2)$ avoiding $\tau$ and containing $\beta$.
The following proposition it the base of all the other
results in this Section, which is obtained by use the same proof argument of 
\cite[Th. 1]{MV3}.

\begin{theorem} 
\label{gencase}
Let $\phi\in S_k(1\mn3\mn2)$ such that $\phi=(\phi^0,m_0,\phi^1,m_1,\cdots,\phi^r,m_r)$, and let 
$\tau=(\tau^0\mn m_0\mn\tau^1\mn m_1\mn\dots\mn\tau^r\mn m_r)$ where $\phi^i=\overline{\tau^i}$. Then 
	$$F_\tau(x)=1+x\sum_{j=0}^r\bigl(F_{\pi^j}(x)-F_{\pi^{j-1}}(x)\bigr)F_{\sigma^j}(x).$$
\end{theorem}

\begin{example} {\rm (Cleasson and Mansour \cite{CM})} 
\label{ex1}
Let us find $F_{12\mn3}(x)$. The canonical decomposition for $12\mn3$ 
is $(12)\mn3$, so $r=0$, and hence $F_{12\mn3}(x)=1+xF_{12}(x)F_{12\mn3}(x)$. 
By definitions, $F_{12}(x)=\frac{1}{1-x}$, hence 
$F_{12\mn3}(x)=\frac{1-x}{1-2x}$. similarly, since $F_{21}(x)=\frac{1}{1-x}$ we 
get that $F_{21\mn3}(x)=\frac{1-x}{1-2x}$.
\end{example}

The case of varying $k$ is more interesting. As an extension of Example \ref{ex1}.

\begin{theorem}
\label{ex11}
Let $\tau\in\{12,1\mn2,21,2\mn1\}$; for any $k\geq 2$, 
 	 $$F_{\tau\mn3\mn4\mn\cdots\mn k}(x)=R_k(x).$$
\end{theorem} 
\begin{proof}
By Theorem \ref{gencase} it is easy to see
	$$F_{\tau\mn3\mn\cdots\mn k}(x)=\frac{1}{1-xF_{\tau\mn3\mn\dots\mn(k-1)}(x)}.$$
Besides, by definitions $F_{\tau}(x)=\frac{1}{1-x}$, hence 
by use \cite[Eq. 4.1]{MV3} the theorem holds.
\end{proof}

\begin{notation} 
The generalized pattern $123\dots k$ we denote by $[k]$.
\end{notation}

\begin{theorem}
\label{avoid1} 
The generating function $F_{[k]}(x)$ is defined by the following equation
	 $$F_{[k]}(x)=\sum_{j=0}^{k-1} (xF_{[k]}(x))^j.$$ 
\end{theorem}
\begin{proof}
Let $\alpha\in S_n(1\mn3\mn2,[k])$; If $\alpha$ nonempty, then we may write 
$\alpha=\beta^1 n\gamma$ such that $\alpha_a=n$ where $\beta^1$ and $\gamma$ 
avoids $1\mn3\mn2$, $\beta^1$ is a permutation of $\{n-1,n-2,\dots,n-a+1\}$, and 
$\gamma$ is a permutation of $\{1,2,\dots,n-a\}$. If $\beta^1$ is nonempty, then 
let us assume that $\beta^1_{a-1}=m_1$ and since $\beta^1$ 
avoids $1\mn3\mn2$, we may write $\beta^1=\delta^1\beta^2 m_1$ where 
$\beta^2$ and $\delta^1$ avoids $1\mn3\mn2$, $\delta^1$ is a permutation 
of $\{n-1,\dots,m_1+1\}$, and $\beta^2$ is a permutation of 
$\{m_1-1,\dots,n-a+1\}$, and so on $\beta^2$. Thus we may 
write $\beta^1=\delta^1\delta^2\dots\delta^{k-2},m_{k-2},\dots,m_1)$ 
where $\delta^j$ avoids $1\mn3\mn2$ and $[k]$ for any $j\leq k-2$, 
and $\delta^{k-1}=\emptyset$.
Therefore, in terms of generating function it is easy to get that 
	$$F_{[k]}(x)=1+xF_{[k]}(x)+\dots+x^{k-1}F_{[k]}^{k-1}(x).$$ 
\end{proof}

\begin{example} {\rm (see Barcucci, Del Lungo, Pergola, and Pinzani \cite{BDPP})}
\label{exx1}
Theorem \ref{avoid1} yields $F_{[1]}(x)=1$, and $F_{[2]}(x)=\frac{1}{1-x}$. The first 
interesting case which Theorem \ref{avoid1} yields $F_{123}(x)=M(x)$. Also
by use Theorem \ref{gencase} we have 
	$$F_{123-4}(x)=\frac{1}{1-xF_{123}(x)}=\frac{1}{2} \sqrt{\frac{1+x}{1-3x}},$$
which is the generating function of Directed animals of size $n$. 
\end{example}

Similarly, by proof of Theorem \ref{avoid1} we get the following.

\begin{theorem}
\label{avoidk}
Let $k>d\geq 1$ and $\tau=\tau'\mn(d+1)(d+2)\dots k$. Then 
    $$F_\tau(x)=\sum_{j=0}^{k-d-1} (xF_\tau(x))^j +x^{k-d}F_\tau^{k-d}(x)F_{\tau'}(x).$$
\end{theorem}

\begin{example} {\rm (Claesson \cite[Pro. 16]{C})}
\label{ex2c}
By Theorem \ref{avoidk} we obtained that $F_{1\mn23}(x)=1+xF_{1\mn23}(x)+x^2F_{1\mn23}(x)$, hence 
$F_{1\mn23}(x)=M(x)$, the generating function of Motzkin numbers. 
\end{example}

The case of varying $k$ is more interesting. As an extension of Example \ref{ex2c}.

\begin{theorem}
\label{ex2c1}
For all $k\geq 3$, the generating functions $F_{1\mn2\mn\cdots\mn(k-2)\mn(k-1)k}(x)$,  
$F_{12\mn3\mn\cdots(k-2)\mn(k-1)k}(x)$, $F_{2\mn1\mn3\mn\cdots\mn(k-2)\mn(k-1)k}(x)$,  
$F_{21\mn3\mn\cdots(k-2)\mn(k-1)k}(x)$ are given by 
	$$\frac{1-x-\sqrt{1-2x+x^2-4x^2R_{k-2}(x)}}{2x^2R_{k-2}(x)}.$$
\end{theorem}
\begin{proof}
By Theorem \ref{avoidk} 
$$\begin{array}{ll}
F_{1\mn2\mn\cdots\mn(k-2)\mn(k-1)k}(x)=&1+xF_{1\mn2\mn\cdots\mn(k-2)\mn(k-1)k}(x)+\\
				       &+x^2F_{1\mn2\mn\cdots\mn(k-2)\mn(k-1)k}^2(x)F_{1\mn2\mn\cdots\mn (k-2)}(x).
\end{array}$$
Therefore, by use Theorem \ref{ex11} the theorem holds in this case. Similarly for    
$12\mn3\mn\cdots\mn(k-2)\mn(k-1)k$, for $2\mn1\mn3\mn\cdots\mn(k-2)\mn(k-1)k$, and 
for $21\mn3\mn\cdots\mn(k-2)\mn(k-1)k$.
\end{proof}

\begin{example}
Theorem \ref{ex2c1} yields that, the generating functions $F_{12\mn34}(x)$, 
$F_{21\mn34}(x)$, $F_{1\mn2\mn34}(x)$, and $F_{2\mn1\mn34}(x)$ are given by 
	$$\frac{1-2x+x^2-\sqrt{1-4x+2x^2+x^4}}{2x^2}.$$
\end{example}

\begin{notation}
The generalized pattern $k\dots 321$ we denote by $<k>$. 
\end{notation}

So, by the same way of proof Theorem \ref{avoid1} we get the following. 

\begin{theorem}
\label{avoid2}
The generating function $F_{<k>}(x)$ is defined by the following equation
	 $$F_{<k>}(x)=\sum_{j=0}^{k-1} (xF_{<k>}(x))^j.$$ 
\end{theorem}
\begin{proof}
Let $\alpha\in S_n(1\mn3\mn2,<k>)$; If $\alpha$ nonempty, then we may write 
$\alpha=\gamma n\beta^1$ such that $\alpha_a=n$ where $\beta^1$ and $\gamma$ 
avoids $1\mn3\mn2$, $\gamma$ is a permutation of $\{n-1,n-2,\dots,n-a+1\}$, and 
$\beta^1$ is a permutation of $\{1,2,\dots,n-a\}$. If $\beta^1$ is nonempty, then 
let us assume that $\beta^1_{1}=m_1$ and since $\beta^1$ 
avoids $1\mn3\mn2$, we may write $\beta^1=m_1 \beta^2\delta^1$ where 
$\beta^2$ and $\delta^1$ avoids $1\mn3\mn2$, $\delta^1$ is a permutation 
of $\{n-a,\dots,m_1+1\}$, and $\beta^2$ is a permutation of 
$\{m_1-1,\dots,1\}$, and so on $\beta^2$. Thus we may 
write $\beta^1=m_1,\dots,m_{k-2},\delta^{k-1},\delta^{k-2}\delta^{k-3}\dots\delta^1)$ 
where $\delta^j$ avoids $1\mn3\mn2$ and $<k>$ for any $j\leq k-2$, 
and since $\alpha$ avoids $<k>$, then $\delta^{k-1}=\emptyset$.
Therefore, in terms of generating function it is easy to get that 
	$$F_{[k]}(x)=1+xF_{[k]}(x)+\dots+x^{k-1}F_{[k]}^{k-1}(x).$$ 
\end{proof}

\begin{example}
\label{exx2}
By Theorem \ref{avoid2} we get $F_{21}(x)=\frac{1}{1-x}$. The first interesting 
special case for $k=3$, 
	$$F_{321}(x)=1+xF_{321}(x)+x^2F_{321}^2(x),$$
hence $F_{321}(x)=M(x)$. Besides, by use Theorem \ref{gencase} we have 
	$$F_{321\mn 4}(x)=\frac{1}{1-xF_{321}(x)}=\frac{1}{2} \sqrt{\frac{1+x}{1-3x}},$$
which is the generating function of Directed animals of size $n$ {\rm (see \cite{BDPP})}. 
\end{example}

Another way to generalize these result by looking at patterns of the form 
$(d+1,\dots,k,1,2,\dots,d)\in S_k$.

\begin{example}
By Theorem \ref{gencase} we get
	$$F_{45\mn 6\mn 12\mn 3}(x)=1+xF_{21}(x)F_{45-6-12-3}(x)+x(F_{45-6-12-3}(x)-F_{12}(x))F_{12-3}(x),$$
which means that 
	$$F_{45\mn 6\mn 12\mn 3}(x)=\frac{(1-3x)(1-x)}{1-5x+6x^2-x^3}.$$
\end{example}

This example can be generalized by use Theorem \ref{gencase} as following.
\begin{theorem}
\label{twol}
Let $k>d$, and let $\tau=\tau'\mn k\mn \tau''\mn d$ be generalized pattern 
such that the minimal letter of $\tau'$ is $d+1$. Then 
	$$F_\tau(x)=\frac{1}{1-x\frac{1-xF_{\tau'}(x)F_{\tau''}(x)}{1-x(F_{\tau'}(x)+F_{\tau''}(x)}}.$$ 
\end{theorem}

\begin{example}
Let $\tau=(d+1)(d+2)\mn(d+3)\mn\dots\mn k\mn12\mn3\mn\dots\mn d$ where $d\geq 2$, $k\geq 4$. By Theorem 
\ref{twol} with use \ref{ex11} 
	$$F_\tau(x)=\frac{1}{1-x\frac{1-xR_{k-d}(x)R_d(x)}{1-x(R_{k-d}(x)+R_d(x)}},$$
and \cite[Section 4]{MV3} yields $F_\tau(x)=R_k(x)$.
\end{example}

For a further generalization of the above example, consider the 
following definition. We say that $\tau$ is a {\it generalized wedge\/} pattern 
if it can be represented as $\tau=(\phi^1,\psi^1,\dots,\phi^r,\psi^r)$ so that 
$\overline{\tau}$ is a {\em wedge} pattern (see \cite{MV3,MV4})
and each $\phi^j$ and $\psi^j$ is either $a\mn(a+1)\mn\dots\mn b$ or $a(a+1)\mn (a+2)\mn\dots\mn b$.
For example, $6\mn 4\mn 5\mn 7\mn 8\mn 3\mn 9\mn 12$ and $45\mn 6\mn 3\mn 7\mn 8\mn 12\mn 9$ are generalized wedge patterns. 
The following theorem is analogue of Mansour and Vainshtein result \cite[Th. 2.6]{MV2}

\begin{theorem} \label{wed}
$F_\tau(x)=R_k(x)$ for any generalized wedge pattern $\tau$ of $k$-letters.
\end{theorem}

\begin{remark}
A comparison of Theorem \ref{wed} with the Main result of \cite{MV2} suggests that 
there should exist a bijection between the sets $S_n(3\mn2\mn1,(m+1)\mn\cdots\mn k\mn1\mn\cdots\mn m)$ 
and $S_n(1\mn3\mn2,\tau)$ for any generalized wedge pattern $\tau$. However, we failed to 
produce such a bijection, and finding it remains a challenging open question.
\end{remark}

\begin{corollary}
\label{ca1}
\begin{enumerate}
\item	$F_{123}(x)=F_{321}(x)=\frac{1-x-\sqrt{1-2x-3x^2}}{2x^2}$;
\item	$F_{132}(x)=\frac{1-\sqrt{1-4x}}{2x}$;
\item	$F_{213}(x)=F_{312}(x)=\frac{1-x^2-\sqrt{(1+x^2)^2-4x}}{2x(1-x)}$;
\item	$F_{231}(x)=\frac{1-x}{1-2x}$.
\end{enumerate}
\end{corollary}
\begin{proof}
The case $(1)$ holds by Theorem \ref{avoid1}, and Theorem \ref{avoid2}; 
the case $(2)$ by definition and Knuths results \cite{Kn}.

To verify $(3)$; let $\alpha\in S_n(213,1\mn3\mn2)$; we may write $\alpha=\beta n \gamma$ such 
that $\alpha_k=n$, where $\beta$ and $\gamma$ are avoids 1\mn3\mn2, $\beta$ is a permutation 
of $\{n-k+1,\dots,n-2,n-1\}$, and $\gamma$ is a permutation of $\{1,2,\dots,n-k\}$. 
Let $m$ the maximal such that the last $m$ numbers in $\beta$ are increasing. So, 
since $\beta$ avoids 1\mn3\mn2 we obtain 
$\beta=(\delta^1,\dots,\delta^m,a^m,\dots,a^1)$ where $a^{j+1}<a^j$ and 
every element in $\delta^j$ is bigger that every element in $\delta^{j+1}$. 
On the other hand, $(\beta,n)$ avoids $213$, so $\delta^m=\emptyset$. Hence 
	$$F_{213}(x)=1+xF_{213}(x)+x^2F_{213}(x)+x^3F_{213}^2(x)+x^4F_{213}^3(x)+\dots,$$
which means that 
	$$F_{213}(x)=1+xF_{213}(x)+\frac{x^2F_{213}(x)}{1-xf_{213}(x)}.$$
Similarly, we obtain $F_{312}(x)=F_{213}(x)$. 

Finally, by use the same argument above proof we easy to get 
$F_{231}(x)=1+x(F_{231}(x)-1)+xF_{231}(x)$, which means that $F_{231}(x)=\frac{1-x}{1-2x}$.
\end{proof}	
\section{Avoiding $1\mn3\mn2$ and containing once a generalized pattern}

Let $g_\tau(n)$ denote the number of permutations in $S_n(1\mn3\mn2)$ that 
contains $\tau$ exactly once where $\tau$ generalized pattern avoids $1\mn3\mn2$, 
and $g_\tau^\phi(n)$ denote the number of permutations in $S_n(1\mn3\mn2,\phi)$ that 
contains $\tau$ exactly once. We denote by $G_\tau(x)$ and
$G_\tau^\phi(x)$ the corresponding ordinary generating functions.

The following statement is the analog of \cite[Th. 3.1]{MV3}, which is 
yields immediately by the same argument proof of \cite[Th. 3.1]{MV3}.

\begin{theorem} 
\label{contain1}
Let $\phi\in S_k(1\mn3\mn2)$ such that the canonical decomposition of $\phi$ 
given by $\phi=(\phi^0,m_0,\phi^1,m_1,\dots,\phi^r,m_r)$, and let 
$\tau=(\theta^0\mn m_0\mn\theta^1\mn m_1\mn\dots\mn\theta^r\mn m_r)$ where $\phi^i=\overline{\theta^i}$. Then 
$$
\bigl(1-xF_{\pi^0}(x)-xF_{\sigma^r}(x)\bigr)G_\tau(x)=
x\sum_{j=1}^r G_{\pi^{j-1}}^{\pi^j}(x)G_{\sigma^j}^{\sigma^{j-1}}(x)
$$
for $r\geq 1$, and
$$
G_\tau(x)=\frac{xF_\tau(x)G_{\pi^0}(x)}{1-xF_{\pi^0}(x)}
$$
for $r=0$.
\end{theorem}

\begin{theorem}
\label{con11}
Let $k\geq 1$,
	$$G_{[k]}(x)=\sum_{j=1}^{k-1} jx^jG_{[k]}(x)F_{[k]}^{j-1}(x) + x^kF_{[k]}^k(x),$$
where $F_{[k]}(x)$ given in Theorem \ref{avoid1}.
\end{theorem}
\begin{proof}
By use the argument proof of Theorem \ref{avoid1} and 
Theorem \ref{avoid2}, the rest is easy to see. 
\end{proof}

\begin{example}
\label{exc1}
Theorem \ref{con11}, and Theorem \ref{ex11} yields $G_{12}(x)=G_{21}(x)=\frac{x^2}{(1-x)^3}$. 
Theorem \ref{con11}, Example \ref{exx1} and Example \ref{exx2} yields 
	$$G_{123}(x)=G_{321}(x)=\frac{x^3M^3(x)}{1-x-2x^2M(x)}=2x-1+\frac{1-3x+2x^3}{2x^3\sqrt{1-2x-3x^2}}.$$ 
\end{example}

\begin{corollary}
\label{cd2}
For any $k\geq 2$; 
   $$G_{12\mn3\mn \cdots\mn k}(x)=\frac{1}{(1-x)U_k^2\left( \frac{1}{2\sqrt{x}} \right) }.$$
\end{corollary}
\begin{proof}
Let $\tau=12\mn3\mn\cdots\mn k$; then $r=0$, and it follows from Theorem \ref{contain1} that
	$$G_{12\mn3\mn\dots\mn k}(x)=\frac{xF_{12\mn3\mn\cdots\mn k}(x)G_{12\mn3\mn\cdots\mn(k-1)}(x)}{1-xF_{12\mn3\mn\dots\mn(k-1)}(x)}.$$
Since $F_{12\mn3\mn\cdots\mn k}(x)=R_k(x)$ (see \ref{ex11}) and $R_k(x)(1-xR_{k-1}(x))=1$, 
we get $G_{12\mn3\mn\dots\mn k}(x)=xR_k^2(x)G_{12\mn3\mn \dots\mn (k-1)}(x)$, which together with 
$G_{12}(x)=x^2R_2^3$ (see \ref{exc1}) gives first part of the theorem. 
\end{proof}

More generally we can get a formula for $G_{12\dots d\mn (d+1)\mn \dots\mn k}(x)$ as follows.

\begin{theorem}\label{gdd1}
Let $k\geq d\geq 1$; 
	$$G_{12\dots d\mn(d+1)\mn(d+2)\mn \dots\mn k}(x)
		=\frac{U_d^2\left( \frac{1}{2\sqrt{x}} \right) }{U_k^2\left( \frac{1}{2\sqrt{x}} \right) }G_{[d]}(x).$$
\end{theorem}
\begin{proof}
By use Theorem \ref{contain1} with induction on $k$ (see Corollary \ref{cd2}) 
it follows that 
	$$G_{12\dots d\mn(d+1)\mn\dots\mn k}(x)=G_{12\dots d}(x)\prod_{j=d+1}^k (xR_j^2(x)).$$
\end{proof}

\begin{example}
\label{ex21x}
Let $H(x)$ be the generating function for number of all $1\mn3\mn2$-avoiding permutations
such that avoiding $21\mn3$ and containing $21$ exactly once; so it is easy to see 
$H(x)=x^2F_{21}(x)(F_{21}(x)-1)$. On the other hand, by Theorem \ref{contain1} we get that 
	$G_{21\mn3}(x)=xH(x)F_{21}(x)+xF_{21}(x)G_{21\mn3}(x)$. 
Hence 
	$$G_{21\mn3}(x)=\frac{1}{U_3^2\left( \frac{1}{2\sqrt{x}} \right)}.$$
\end{example}

This example can be generalized by use Theorem \ref{contain1} as following.

\begin{theorem}
\label{g21}
For any $k\geq 3$,
	$$G_{21\mn3\mn\cdots\mn k}(x)=\frac{1}{U_k^2\left( \frac{1}{2\sqrt{x}} \right)}.$$
\end{theorem}
\begin{proof}
By Theorem \ref{contain1} we get that for $k\geq 4$
	$$G_{21\mn3\mn\cdots\mn k}(x)=\frac{xG_{21\mn3\mn\cdots\mn (k-1)}(x)F_{21\mn3\mn\cdots\mn k}(x)}{1-xF_{21\mn3\mn\cdots\mn (k-1)}(x)},$$
and by use Theorem \ref{ex11} we get that 
	$$G_{21\mn3\mn\cdots\mn k}(x)=\frac{xR_k(x)}{1-xR_{k-1}(x)}G_{21\mn3\mn\cdots\mn (k-1)}(x),$$
hence, by use the identity $R_k(x)(1-xR_{k-1}(x))=1$ we have that 
	$$G_{21\mn3\mn\cdots\mn k}(x)=xR_k^2(x)G_{21\mn3\mn\cdots\mn (k-1)}(x).$$ 
By use induction on $k$ and Example \ref{ex21x} the theorem holds.
\end{proof}

One can try to obtain results similar to Theorems \ref{cd2}, \ref{gdd1} and \ref{g21}, but 
expressions involved become extremely cumbersome. So we just consider a simplest cases.
\section{Containing $1\mn3\mn2$ exactly once and avoiding another generalized pattern}
Denote by $H_{\tau}(x)$ the generating function for the number
of permutations in $S_n$ avoiding a generalized pattern $\tau$ 
and containing $1\mn3\mn2$ exactly once. 
We start from the following result, which is analog of result 
obtained in \cite{MV1}.

\begin{theorem}
\label{h1}
For any $k\geq 2$,
	$$H_{12\mn3\mn\cdots\mn k}(x)=\frac{x}{U_k^2\left( \frac1{2\sqrt{x}} \right)} \sum_{j=1}^{k-2} U_j^2\left( \frac1{2\sqrt{x}} \right)$$
\end{theorem}
\begin{proof}
Let $\alpha\in S_n$ contain $132$ exactly once. There are two
possibilities: either the only occurrence of $132$ in $\alpha$ contains $n$, 
or does not contain it. In the first case we get immediately that any entry
of $\alpha$ to the right of $n$ is less than any entry of $\alpha$ to the left
of $n$, since otherwise one gets an occurrence of $132$ involving $n$.
In the second case, let $i\; n\; j$ be the occurrence of $132$ in $\alpha$. 
First,
we have $j=i+1$, since otherwise either $i+1\; n \;j$ or $i\; n\; i+1$ would 
be a
second occurrence of $132$ in $\alpha$. Next, $i$ immediately precedes $n$ 
in $\alpha$, since if $l$ lies between $i$ and $n$, then either $l\; n\; j$ or
$i\;l\;j$ would be a second occurrence of $132$ in $\alpha$. Finally, any entry
to the right of $j$ is less than any entry between $n$ and $j$, which in its 
turn, is less than any entry to the left of $i$ (the proof is similar to the
analysis of the first case). Therefore, we have exactly three possibilities for
the structure of $\alpha$:

(i) there exists $t$ such that $\alpha=(\alpha',n,\alpha'')$, where 
$\alpha'$ is a permutation of $n-1,n-2,\dots,n-t+1$ containing
$132$ exactly once, and $\alpha''$ is a permutation of  
$1,2,\dots,n-t$ avoiding $132$;

(ii) there exists $t$ such that $\alpha=(\alpha',n,\alpha'')$, where 
$\alpha'$ is a permutation of $n-1,n-2,\dots,n-t+1$ avoiding
$132$, and $\alpha''$ is a permutation of  
$1,2,\dots,n-t$ containing $132$ exactly once;

(iii) there exist $t,u$ such that 
$\alpha=(\alpha',n-t+1,n,\alpha'',n-t+2,\alpha''')$, where $\alpha'$ 
is a permutation of $n-1,n-2,\dots,n-t+3$ avoiding $132$, $\alpha''$ is a
permutation of $n-t,n-t-1,\dots,n-u+1$ avoiding $132$, and $\alpha'''$ is a 
permutation of $1,2,\dots,n-u$ avoiding $132$.

We are ready to write an equation for $H_{12\mn3\mn\cdots\mn k}(x)$ with $k\geq 3$.
The contribution of the first structure above is
$xH_{12\mn3\mn\cdots\mn (k-1)}(x) F_{12\mn3\mn\cdots\mn k}(x)$. 
The contribution of the second possible structure is 
$x\bar F_{12\mn3\mn\cdots\mn(k-1)}(x)H_{12\mn3\mn\cdots\mn k}(x)$.
Finally, the contribution of the third possible structure is
$x^3\bar F_{12\mn3\mn\cdots (k-1)}^2(x)\bar F_{12\mn3\mn\cdots\mn k}(x)$.

Solving the obtained linear equation and using Theorem \ref{ex11}, $H_{12}(x)=0$ 
and well known identities involving Chebyshev polynomials (see e.g. \cite{MV2}), 
we get the desired expression for $H_{12\mn3\mn\cdots\mn k}(x)$, $k\geq 3$.
\end{proof}

\begin{example}
Theorem \ref{h1} yields $H_{12\mn3}(x)=\frac{1}{U_3^2\left( \frac1{2\sqrt{x}} \right)}$. 
Similarly by the structure of permutations containing exactly $1\mn3\mn2$ {\rm (see proof of Theorem \ref{h1})}  
we obtain that $H_{21\mn3}(x)=\frac{1-x}{U_3^2\left( \frac1{2\sqrt{x}} \right)}$.
\end{example}

By this example with use the proof of Theorem \ref{h1}, we can be extended this example 
to generalized pattern $21\mn3\mn\cdots\mn k$

\begin{theorem}
\label{h21}
For all $k\geq 3$,
	$$H_{21\mn3\mn\cdots\mn k}(x)=\frac{x}{U_k^2\left( \frac1{2\sqrt{x}} \right)} \left( \sum_{j=1}^{k-2} U_j^2\left( \frac1{2\sqrt{x}} \right) -1\right)$$
\end{theorem}

One can try to obtain results similar to Theorems \ref{h1} and \ref{h21}, but 
expressions involved become extremely cumbersome. So we just consider a two 
simplest cases.
\section{Containing $1\mn3\mn2$ and another generalized pattern exactly once}
Denote by $\Phi_{\tau}(x)$ the generating function for the number
of permutations in $\SS_n$ containing both $1\mn3\mn2$ and a generalized 
pattern $\tau$ exactly once. We start from the following result.

\begin{theorem} For any $k\geq 2$,
$$
\Phi_{12\mn3\mn\cdots\mn k}(x)=\frac{1}{U_2(t)U_k^2(t)}\left[ 1+\sum_{i=2}^{k-1} 
\frac{2\sqrt{x}}{U_i(t)U_{i+1}(t)} \left( \sum_{j=1}^i U_j^2(t) -1 \right)\right], \qquad t=\frac1{2\sqrt{x}}.
$$
\end{theorem}
\begin{proof} Let $\beta^k=\tau\mn3\mn\cdots\mn k$. 
The three possible structures of permutations containing $1\mn3\mn2$
exactly once are described in the proof of Theorem \ref{h1}. Let us find the 
recursion for $\Phi_{\beta^k}(x)$. It is easy to see
that the contribution of the first structure equals
$$x\Phi_{\beta^{k-1}}(x) F_{\beta^k}(x)+xH_{\beta^{k-1}}(x)G_{\beta^k}(x),$$
the contribution of the second structure equals
$$xG_{\beta^{k-1}}(x)H_{\beta^k}(x)+xF_{\beta^{k-1}}(x)\Phi_{\beta^k}(x),$$
while the contribution of the third structure equals
$$2x^3G_{\beta^{k-1}}(x) F_{\beta^{k-1}}(x) F_{\beta^k}(x)+x^3 F_{\beta^{k-1}}^2(x)G_{\beta^k}(x).$$ 
Solving the obtained recursion with the initial condition $\Phi_{12}(x)=\frac{x^3}{(1-x)^3}$ 
(which easy to compute by definitions), and using Theorems \ref{h1}, \ref{ex11}, and \ref{cd2},
we get the desired result.
\end{proof}

Another interesting results is $\Phi_{21\mn3\mn\cdots\mn k}(x)$. Similarly to 
proof of the above theorem with initial condition $\Phi_{21\mn3}(x)=\frac{2x^4(1-x)^2}{(1-2x)^3}$
(as a remark $\Phi_{21}(x)=\frac{x^3}{(1-x)^2}$), and using 
Theorems \ref{ex11}, \ref{g21} and \ref{h21}, we obtain the following.

\begin{theorem}
$$
\Phi_{21\mn3\mn\cdots\mn k}(x)=\frac{1}{4t^3 U_k^2(t)}\left[ \frac{U_2^2(t)}{U_3(t)} +\sum_{i=3}^{k-1} 
\frac{1}{U_i(t)U_{i+1}(t)} \left( \sum_{j=1}^i U_j^2(t) -2 \right) \right], \qquad t=\frac1{2\sqrt{x}}.
$$
\end{theorem}

\begin{remark}
In the above sections, we calculated either $F_\tau(x)$, or $G_\tau(x)$, or $H_\tau(x)$, or 
$\Phi_\tau(x)$ for given $\tau$ generalized patterns of $k$-letters. We note that there 
are others interesting cases of $\tau$ which we can calculate these functions in 
the same methods. So we just consider a simplest and interesting cases.
\end{remark}


\begin{thebibliography}{99}
\bibitem[B]{B1}
{\sc M. B{\"o}na}, The permutation classes equinumerous to the smooth class, {\em 
Electron. J. Combin.} {\bf 5} (1998), \#R31

\bibitem[BS]{BS} 
{\sc E. Babson and E. Steingrimsson}, Generalized permutation patterns and 
a classification of the Mahonian statistics, {\em S\'eminaire Lotharingien de 
Combinatoire}, B44b:18pp, (2000).

\bibitem[BDPP]{BDPP}
{\sc E, Barcucci, A. Del Lungo, E. Pergola, and R. Pinzani}, 
Directed animals, forests and permutations, 
{\em Disc. Math.} {\bf 204} (1999) 41--71.

\bibitem[C]{C}
{\sc A. Claesson}, Generalised pattern avoidance, 
{\em European Journal of Combinatorics}, {\bf 22} (2001) 961--973.

\bibitem[CM]{CM}
{\sc A. Claesson and T. Mansour}, Permutations avoiding a pair of generalized patterns of 
length three with exactly one dash, preprint CO/0107044.

\bibitem[CW]{CW}    
{\sc T. Chow and J. West}
Forbidden subsequences and Chebyshev polynomials
{\em Discr. Math.} {\bf 204} (1999) 119--128.

\bibitem[Ki]{Ki}
{\sc S. Kitaev}, Multi-Avoidance of generalised patterns, to appear.

\bibitem[Kn]{Kn}
{\sc D.E. Knuth}, {\em The Art of Computer Programming}, 2nd ed. Addison 
Wesley, Reading, MA (1973).

\bibitem[Kr]{Kr}
{\sc C. Krattenthaler}, Permutations with restricted patterns and 
Dyck paths (2000), preprint CO/0002200

\bibitem[K]{Km}
{\sc D. Kremer}, Permutations with forbidden subsequnces and a generalized 
Schr\"oder number, {\em Disc. Math.} {\bf 218} (2000), 121--130.

\bibitem[MV1]{MV1}
{\sc T. Mansour and A. Vainshtein}, 
Restricted permutations, continued fractions, and Chebyshev polynomials, 
{\em The Electronic Journal  of Combinatorics} {\bf 7} (2000), \#R17.

\bibitem[MV2]{MV3}
{\sc T. Mansour and A. Vainshtein}, 
Restricted 132-avoiding permutations, {\em Adv. Appl. Math.}  {\bf 126} (2001), 258--269.

\bibitem[MV3]{MV2}
{\sc T. Mansour and A. Vainshtein}, 
Layered restrictions and Chebychev polynomials (2000), 
{\em Annals of Combinatorics}, to appear (2001), preprint CO/0008173.

\bibitem[SS]{SS}
{\sc R. Simion and F.W. Schmidt,} Restricted permutations, 
{\em European Journal of Combinatorics} {\bfseries 6} (1985) 383--406.

\bibitem[R]{R} 
{\sc A. Robertson}, Permutations containing and avoiding $123$ and $132$ patterns, 
{\em Disc. Math. and Theo. Comp. Sci.} {\bf 3} (1999), 151--154.

\bibitem[RWZ]{RWZ}
{\sc A. Robertson, H. Wilf, and D. Zeilberger}, Permutation patterns and continuous 
fractions, {Electron. J. Combin.} {\bf 6} (1999), \#R38.

\bibitem[W]{W}
{\sc J. West}, Generating trees and forbidden subsequences, {\em Disc. Math.} {\bf 157} 
(1996), 363--372.
\end{thebibliography}
\end{document}